\documentclass[10pt]{amsart}

\usepackage{amsmath,amsthm,amsfonts,amscd,amssymb}
\usepackage{hyperref,wasysym}
\usepackage{verbatim}
\usepackage{amssymb}

\newtheorem{thm}{Theorem}[section]
\newtheorem{cor}[thm]{Corollary}
\newtheorem{lem}[thm]{Lemma}

\theoremstyle{definition}

\theoremstyle{remark}

\begin{document}

\title{Integral zeros of quadratic polynomials avoiding sublattices}

\author{Lenny Fukshansky}
\author{Sehun Jeong}

\address{Department of Mathematics, 850 Columbia Avenue, Claremont McKenna College, Claremont, CA 91711}
\email{lenny@cmc.edu}
\address{Institute of Mathematical Sciences, Claremont Graduate University, Claremont, CA 91711}
\email{sehun.jeong@cgu.edu}


\begin{abstract} 
Assuming an integral quadratic polynomial with nonsingular quadratic part has a nontrivial zero on an integer lattice outside of a union of finite-index sublattices, we prove that there exists such a zero of bounded norm and provide an explicit bound. This is a contribution related to the celebrated theorem of Cassels on small-height zeros of quadratic forms, which builds on some previous work in this area. We also demonstrate an application of these results to the problem of effective distribution of angles between vectors in the integer lattice.
\end{abstract}

\maketitle

\def\A{{\mathcal A}}
\def\B{{\mathcal B}}
\def\C{{\mathcal C}}
\def\D{{\mathcal D}}
\def\F{{\mathcal F}}
\def\x{{\mathcal H}}
\def\I{{\mathcal I}}
\def\J{{\mathcal J}}
\def\K{{\mathcal K}}
\def\L{{\mathcal L}}
\def\M{{\mathcal M}}
\def\N{{\mathcal N}}
\def\O{{\mathcal O}}
\def\R{{\mathcal R}}
\def\s{{\mathcal S}}
\def\V{{\mathcal V}}
\def\W{{\mathcal W}}
\def\X{{\mathcal X}}
\def\Y{{\mathcal Y}}
\def\H{{\mathcal H}}
\def\Z{{\mathcal Z}}
\def\OO{{\mathcal O}}
\def\BB{{\mathbb B}}
\def\cee{{\mathbb C}}
\def\EE{{\mathbb E}}
\def\Nn{{\mathbb N}}
\def\pee{{\mathbb P}}
\def\que{{\mathbb Q}}
\def\real{{\mathbb R}}
\def\zed{{\mathbb Z}}
\def\hyp{{\mathbb H}}
\def\aa{{\mathfrak a}}
\def\HH{{\mathfrak H}}
\def\qbar{{\overline{\mathbb Q}}}
\def\eps{{\varepsilon}}
\def\ahat{{\hat \alpha}}
\def\bhat{{\hat \beta}}
\def\gt{{\tilde \gamma}}
\def\h{{\tfrac12}}
\def\be{{\boldsymbol e}}
\def\bei{{\boldsymbol e_i}}
\def\bff{{\boldsymbol f}}
\def\ba{{\boldsymbol a}}
\def\bb{{\boldsymbol b}}
\def\bc{{\boldsymbol c}}
\def\bm{{\boldsymbol m}}
\def\bk{{\boldsymbol k}}
\def\bi{{\boldsymbol i}}
\def\bl{{\boldsymbol l}}
\def\bq{{\boldsymbol q}}
\def\bu{{\boldsymbol u}}
\def\bt{{\boldsymbol t}}
\def\bs{{\boldsymbol s}}
\def\bv{{\boldsymbol v}}
\def\bw{{\boldsymbol w}}
\def\bx{{\boldsymbol x}}
\def\bX{{\boldsymbol X}}
\def\bz{{\boldsymbol z}}
\def\bwy{{\boldsymbol y}}
\def\bY{{\boldsymbol Y}}
\def\bL{{\boldsymbol L}}
\def\baa{{\boldsymbol\alpha}}
\def\bbb{{\boldsymbol\beta}}
\def\bet{{\boldsymbol\eta}}
\def\bxi{{\boldsymbol\xi}}
\def\bo{{\boldsymbol 0}}
\def\bol{{\boldkey 1}_L}
\def\ep{\varepsilon}
\def\p{\boldsymbol\varphi}
\def\q{\boldsymbol\psi}
\def\rank{\operatorname{rank}}
\def\aut{\operatorname{Aut}}
\def\lcm{\operatorname{lcm}}
\def\sgn{\operatorname{sgn}}
\def\spn{\operatorname{span}}
\def\md{\operatorname{mod}}
\def\Norm{\operatorname{Norm}}
\def\dim{\operatorname{dim}}
\def\det{\operatorname{det}}
\def\Vol{\operatorname{Vol}}
\def\rk{\operatorname{rk}}
\def\Gal{\operatorname{Gal}}
\def\WR{\operatorname{WR}}
\def\WO{\operatorname{WO}}
\def\GL{\operatorname{GL}}
\def\pr{\operatorname{pr}}
\def\Tr{\operatorname{Tr}}

\section{Introduction and statement of results}
\label{intro}

Let $n \geq 2$ and let
$$F(\bx,\bwy) = \sum_{i=1}^n \sum_{j=1}^n f_{ij} x_i y_j \in \zed[\bx,\bwy]$$
be a symmetric bilinear form in $2n$ variables with integer coefficients $f_{ij} = f_{ji}$, and let $F(\bx) := F(\bx,\bx)$ be the corresponding integral quadratic form in $n$ variables. We say that $F$ is {\it isotropic} over $\zed$ if there exists a point $\bo \neq \bz \in \zed^n$ such that $F(\bz)=0$. A classical 1955 theorem of Cassels~\cite{cassels:small} (see also \S 6.8 of~\cite{cassels_qf}) asserts that an isotropic integral quadratic form in $n$ variables has a nontrivial integral zero $\bz$ of small size. Specifically, we will use the norms
$$|\bz| = \max_{1 \leq i \leq n} |z_i|,\ \|\bz\| = \left( z_1^2 + \dots + z_n^2 \right)^{1/2}$$
to measure the size of $\bz$, so that $|\bz| \leq \|\bz\| \leq \sqrt{n}\ |\bz|$. Then the bound obtained by Cassels is of the form
\begin{equation}
\label{cassels_bnd}
|\bz| \ll |F|^{\frac{n-1}{2}},
\end{equation}
where $|F| = \max_{1 \leq i,j \leq n} |f_{ij}|$ and the constant in the Vinogradov notation $\ll$ depends only on $n$. Cassels' theorem has opened a lively new avenue of research into the effective arithmetic theory of quadratic forms; see~\cite{cassels_overview} for a survey of many results by different authors in this general direction.

Notice that the questions of existence of integral or rational zeros for the quadratic form $F$ are equivalent. On the other hand, these questions become quite different for inhomogeneous quadratic equations. Let us write $\F = (f_{ij})_{1 \leq i,j \leq n}$ for the $n \times n$ symmetric coefficient matrix of $F$, then
$$F(\bx,\bwy) = \bx^{\top} \F \bwy.$$
From here on, we will assume that the form $F$ is {\it regular}, meaning that the coefficient matrix $\F$ is nonsingular. Define an inhomogeneous quadratic polynomial in $n$ variables $\bx = (x_1,\dots,x_n)$ as
$$Q(\bx) = F(\bx) + L(\bx) + t,$$
where $F(\bx)$ is a quadratic form as above, $L(\bx) = \sum_{i=1}^n \ell_i x_i$ is a linear form with integer coefficients, and $t \in \zed$. Since $F$ is regular, we will refer to this $Q$ as regular too. We write $|L|$ for $\max_{1 \leq i \leq n} |\ell_i|$ and set $|Q| = \max \{ |F|, |L|, |t| \}$. Masser in~\cite{masser} proved the existence of small-size rational solutions for an equation $Q(\bx) = 0$ with the bound being in terms of $|Q|$, assuming that $Q$ is isotropic over $\que$. Our main interest in this paper, however, is in integer solutions. Assuming that $n \geq 3$ and $Q$ is isotropic over $\zed$, Dietmann proved in~\cite{dietmann} that there exists $\bz \in \zed^n$ such that $Q(\bz) = 0$ and
\begin{equation}
\label{RD_bnd}
|\bz| \ll |Q|^{\rho(n)},
\end{equation}
where
\begin{equation}
\label{pk_bnd}
\rho(n) = \begin{cases}
2100 & \mbox{ if $n = 3$},\\
84 & \mbox{ if $n = 4$},\\
5n + 19 + 74/(n - 4) & \mbox{ if $n \geq 5$}
\end{cases}
\end{equation}
and the constant in upper bound of~\eqref{RD_bnd} depends only on $n$ and is effectively computable. In the case $n=2$, Kornhauser~\cite{kornhauser} proved under the same conditions that $Q(\bx) = 0$ has an integer solution $\bz$ with
\begin{equation}
\label{kornh}
|\bz| \leq (28|Q|)^{10|Q|},
\end{equation}
and showed that in the binary case an upper bound on $|\bz|$ that would be polynomial in $Q$ is, in general, not possible.

We want to focus on the distribution of small-size zeros of integral quadratic polynomials. Specifically, one can speculate that if these zeros are ``well-distributed", in some sense, it should not be easy to ``cut them out" by a finite collection of sublattices of the integer lattice. In particular, Theorem~1.5 of~\cite{lf_chan} implies that if a regular isotropic integral quadratic form $Q$ assumes the value $t \in \zed$ on a rank-$k$ lattice $\Lambda \subseteq \zed^n$, $3 \leq k \leq n$ and $\Omega_1,\dots,\Omega_m \subset \Lambda$ are sublattices of rank $k-1$, then $Q(\bz) = t$ for some small-size point $\bz \in \Lambda \setminus \left( \bigcup_{i=1}^m \Omega_i \right)$, with an explicit bound on $|\bz|$. In other words, the equation $Q(\bx) = t$ has solutions of controllably small size avoiding any finite union of sublattices of smaller rank. A key observation implicitly used in the proof of this result is a certain ``projective nature" of the problem: there exist points $\bz \in \Lambda$ such that $\alpha \bz \not\in \left( \bigcup_{i=1}^m \Omega_i \right)$ for any $\alpha \in \zed$. This is due to the assumption that the sublattices $\Omega_1,\dots,\Omega_m$ have smaller rank than $\Lambda$. On the other hand, if their rank were also $k$, then all of them would have finite index in $\Lambda$ and hence for any $\bz \in \Lambda$ there exist integers $\alpha_1,\dots,\alpha_m$ so that $\alpha_i \bz \in \Omega_i$ for every $1 \leq i \leq m$. This observation makes the method of~\cite{lf_chan} unusable for the case of sublattices of finite index, which constitute an ``inhomogeneous" situation, in a certain sense. The main result of our present paper addresses precisely this situation using a rather different method.

\begin{thm} \label{main} Let $\Lambda \subseteq \zed^n$ be a sublattice of rank $k$, $2 \leq k \leq n$. Let $\Omega_1,\dots,\Omega_m \subset \Lambda$ be proper sublattices of finite indices and let $\Omega = \bigcap_{j=1}^m \Omega_j$ be their intersection sublattice, which then also has a finite index in $\Lambda$. Assume that $Q$ is regular and there exists some point $\bz \in \Lambda \setminus \left( \bigcup_{j=1}^m \Omega_j \right)$ such that $Q(\bz) = 0$. If $k \geq 3$, then there exists such $\bz$ with
\begin{equation}
\label{main_bnd}
|\bz| \ll \det (\Omega)^{2\rho(k)+1} |Q|^{\rho(k)},
\end{equation}
where $\rho(k)$ is given by \eqref{pk_bnd} with $n$ replaced by $k$. If $k=2$, then there exists such $\bz$ with
\begin{equation}
\label{main_bnd-2}
|\bz| \ll \det (\Omega) \left( 2408 (\det (\Omega))^2 |Q| \right)^{860 n^2 (\det (\Omega))^2 |Q|},
\end{equation}
and the constants in upper bounds of~\eqref{main_bnd} and~\eqref{main_bnd-2} depend only on $n$ and $k$ and are effectively computable.
\end{thm}

\noindent
We prove Theorem~\ref{main} in Section~\ref{main_proof}. Our argument uses Hermite's inequality to select a ``short" basis for $\Lambda$ and then a variant of a normal form for a corresponding basis for the intersection lattice $\Omega$. This allows us to choose ``short" coset representatives of $\Omega$ in $\Lambda$ and reduce our problem to a search for small-size integer zero of a regular integral quadratic polynomial. We then use Dietmann's theorem.

Our method also allows for a simple proof of existence of a small-size point in our lattice $\Lambda$ outside of a finite union of finite-index sublattices. A more general version of this problem was previously considered by Henk and Thiel in~\cite{henk_thiel} and their bound is generally better than ours. We still record our observation in Corollary~\ref{one_out} only because it is a particularly simple alternative proof.

In Section~\ref{angle_sec}, we show an application of the bounds on zeros of quadratic equations discussed above to the problem of effective distribution of angles between vectors in the integer lattice~$\zed^n$. Specifically, fixing a vector $\ba \in \zed^n$ and assuming that there exists some vector $\bb \in \zed^n$ making a prescribed angle $\theta$ with $\ba$, one can ask for such a vector of bounded norm. A similar question can be asked with the additional condition that the vector $\bb$ avoids a union of finite-index sublattices. We discuss these questions and derive respective bounds in Corollaries~\ref{angle_eff-1} and~\ref{angle_eff-2}. We are now ready to proceed.
\bigskip

\section{Proof of Theorem~\ref{main}}
\label{main_proof}

Let us set $\Delta = \det(\Lambda)$ and start by selecting a basis for $\Lambda$. Hermite's inequality (see, e.g., \cite{martinet}, Section~2.2) guarantees that there exists a basis $\ba_1,\dots,\ba_k$ for $\Lambda$ such that
\begin{equation}
\label{hermite}
\max_{1 \leq i \leq k} |\ba_i| \leq \prod_{i=1}^k |\ba_i| \leq \prod_{i=1}^k \|\ba_i\| \leq \left( \frac{4}{3} \right)^{\frac{k(k-1)}{2}} \Delta.
\end{equation}
The first inequality follows from the fact that each $\ba_i \in \zed^n$, and hence $|\ba_i| \geq 1$. Now, there exists a basis $\bb_1,\dots,\bb_k$ for $\Omega$ so that
\[ \left\{ \begin{array}{ll}
\bb_1 = v_{11} \ba_1 \\
\bb_2 = v_{21} \ba_1 + v_{22} \ba_2 \\
\dots \dots \dots \dots \dots \dots \dots \dots \\
\bb_k = v_{k1} \ba_1 + \dots + v_{kk} \ba_k,
\end{array} \right. \]
where all $v_{ij} \in \zed$ and $0 \leq v_{ij} < v_{ii}$ for all $1 \leq j < i \leq k$ (see Section~I.2 of~\cite{cassels_gn}).  Let us write $d_i = [\Lambda : \Omega_i] > 1$ for each $1 \leq i \leq m$, then $d := [\Lambda : \Omega] \leq d_1 \cdots d_m$. With this notation, we have 
$$\prod_{i=1} v_{ii} = d = \left| \Lambda / \Omega \right|$$
is the number of cosets of $\Omega$ in $\Lambda$. In particular, since all the $v_{ij}$ are positive integers, this implies that
$$\max_{1 \leq i,j \leq k} v_{ij} \leq d.$$
Let us write $\bc_1,\dots,\bc_d$ for the coset representatives of the form $\bc_1 = \bo$ and 
\begin{equation}
\label{cosets}
\bo \neq \bc_i = \sum_{j=1}^k q_{ij} \ba_j,\ 0 \leq q_{ij} < v_{ii},
\end{equation}
for all $2 \leq i \leq d$. Then notice that
\begin{equation}
\label{bb_bnd}
\max_{1 \leq i \leq k} |\bb_i| \leq k \left( \max_{1 \leq i \leq k} v_{ii} \right) \left( \max_{1 \leq i \leq k} |\ba_i| \right) \leq \left( \frac{4}{3} \right)^{\frac{k(k-1)}{2}}  kd \Delta,
\end{equation}
and analogously
\begin{equation}
\label{cc_bnd}
\max_{2 \leq i \leq k} |\bc_i| \leq k \left( \max_{1 \leq i,j \leq k} q_{ij} \right) \left( \max_{1 \leq i \leq k} |\ba_i| \right) \leq \left( \frac{4}{3} \right)^{\frac{k(k-1)}{2}} kd \Delta.
\end{equation}
Let us write $B = (\bb_1 \dots \bb_k)$ for the $n \times k$ matrix, whose columns are the basis vectors $\bb_1,\dots,\bb_k$ for $\Omega$. Notice that every $\bx \in \Lambda$ can be written in the form
\begin{equation}
\label{x_decomp}
\bx = \bc_i + \sum_{j=1}^k x_j \bb_j = \bc_i + B \bx',
\end{equation}
for some $2 \leq i \leq d$ and $\bx' := (x_1,\dots,x_k)^{\top} \in \zed^k$. Since $B\bx' \in \Omega_j$ for every $1 \leq j \leq m$, we have $\bx \not\in \bigcup_{j=1}^m \Omega_j$ if and only if the corresponding $\bc_i \not\in \bigcup_{j=1}^m \Omega_j$. Hence, the equation
\begin{eqnarray*}
G_i(x_1,\dots,x_k) & := & Q \left(  \bc_i + \sum_{j=1}^k x_j \bb_j \right) \\
& = & \sum_{r=1}^n \sum_{s=1}^n f_{rs} \left( c_{ir} + \sum_{j=1}^k b_{jr} x_j \right) \left( c_{is} + \sum_{j=1}^k b_{js} x_j \right) \\
&\ & + \sum_{r=1}^n \ell_r \left( c_{ir} + \sum_{j=1}^k b_{jr} x_j \right) + t = 0
\end{eqnarray*}
has a solution in integers $x_1,\dots,x_k$ for some $2 \leq i \leq d$. Fix this $i$, and let us crudely estimate $|G_i|$, using~\eqref{bb_bnd} and~\eqref{cc_bnd}:
\begin{eqnarray}
\label{gt_bnd}
|G_i| & \leq & \left( n^2 (k+1)^2 + n(k+1) + 1 \right) |Q| \max \left\{ |\bc_i|^2, |\bc_i| \max_{1 \leq j \leq k} |\bb_j|, \max_{1 \leq j \leq k} |\bb_j|^2 \right\}  \nonumber \\
& \leq & \left( \frac{4}{3} \right)^{\frac{2k(k-1)}{2}}  \left( n^2 (k+1)^2 + n(k+1) + 1 \right)  k^2 d^2 \Delta^2 |Q|.
\end{eqnarray}
Now, $G_i(\bx') = 0$ is an inhomogeneous quadratic equation in $k$ variables with integer coefficients which has an integer solution. Further, the quadratic part of this equation can be written as
$$\bx'^{\top} (B^{\top} \F B) \bx',$$
where the $k \times k$ coefficient matrix $B^{\top} \F B$ is nonsingular. Then, by a theorem of Dietmann~\cite{dietmann} (see inequality~\eqref{RD_bnd} above), there exists a point $\bz' \in \zed^k$ such that $G_i(\bz') = 0$ and 
\begin{equation}
\label{dietmann_bnd}
|\bz'| \ll \begin{cases}
(28|G_i|)^{10|G_i|} & \mbox{ if $k = 2$},\\
|G_i|^{\rho(k)} & \mbox{ if $k \geq 3$},
\end{cases}
\end{equation}
where the implied constant is $1$ if $k=2$ and depends only on $k$ if $k \geq 3$, and $\rho(k)$ as in~\eqref{pk_bnd}. The corresponding $\bz$ as in~\eqref{x_decomp} is then a zero of $Q$ contained in~$\Lambda \setminus \left( \bigcup_{j=1}^m \Omega_j \right)$. Combining~\eqref{x_decomp} with~\eqref{cc_bnd} and~\eqref{bb_bnd}, it follows that
\begin{equation}
\label{x_bnd}
|\bz| \leq |\bc_i| + |B\bz'| \leq \left( \frac{4}{3} \right)^{\frac{k(k-1)}{2}} kd \Delta + k |B| |\bz'| \leq \left( \frac{4}{3} \right)^{\frac{k(k-1)}{2}} kd \Delta (1 + k |\bz'|).
\end{equation}
Observe that $d\Delta = \det(\Omega)$. The bound~\eqref{main_bnd} now follows upon combining~\eqref{x_bnd} with~\eqref{dietmann_bnd} and~\eqref{gt_bnd}. In the case $k=2$, we presented a slightly weaker bound than actually follows from our inequalities in the interest of a simpler looking result. This completes the proof of Theorem~\ref{main}.
$\ \ \ \ \ \ \ \ \ \ \ \ \ \ \ \ \ \ \ \ \ \ \ \ \ \ \ \ \ \ \ \ \ \ \ \ \ \ \ \ \ \ \ \ \ \ \ \ \ \ \ \ \ \ \ \ \ \ \ \ \ \ \ \ \ \ \ \ \ \ \ \ \ \ \ \square$
\bigskip

We also want to remark that our method allows for a simple observation on the basic problem of finding short vectors in a lattice outside of a collection of sublattices. In the case of sublattices of lower rank, this problem was originally treated in~\cite{lf1}. More recently, Henk and Thiel~\cite{henk_thiel} considered this problem in the case of sublattices of finite index. Specifically, Theorem~1.2 of~\cite{henk_thiel} applied to our situation states that, assuming $\Lambda \not\subseteq \bigcup_{i=1}^m \Omega_i$, there exists $\bz \in \Lambda \setminus \bigcup_{i=1}^m \Omega_i$ with
\begin{equation}
\label{ht}
|\bz| < \frac{\det(\Omega)}{\lambda_1(\Omega)^{k-1}} \left( \sum_{i=1}^m \frac{1}{d_i} - \frac{m - 1}{d} + \frac{\lambda_1(\Omega)^k}{\det(\Omega)} \right),
\end{equation}
where $\lambda_1(\Omega) = \min \{ |\bx| : \bx \in \Omega \setminus \{\bo\} \}$ is the first successive minimum of $\Omega$ with respect to the sup-norm $|\ |$. This result was obtained using a careful analysis and volume computations on the torus group $\real^k/\Omega$. On the other hand, our proof of Theorem~\ref{main} suggests a very simple argument producing a bound for such a point~$\bz$, albeit with weaker than~\eqref{ht}.

\begin{cor} \label{one_out} Let $\Lambda \subseteq \zed^n$ be a sublattice of rank $k$, $2 \leq k \leq n$, and let $\Omega_1,\dots,\Omega_m \subset \Lambda$ be sublattices of finite indices. Let $\Omega = \bigcap_{j=1}^m \Omega_j$, so $\Omega \subset \Lambda$ is also a sublattice of finite index. Assume that $\Lambda \not\subseteq \bigcup_{j=1}^m \Omega_j$, then there exists a point $\bz \in \Lambda \setminus  \left( \bigcup_{j=1}^m \Omega_j \right)$ such that
$$|\bz| \leq \left( \frac{4}{3} \right)^{\frac{k(k-1)}{2}} k \det(\Omega).$$
\end{cor}

\proof
Since $\Lambda \not\subseteq \bigcup_{j=1}^m \Omega_j$, at least one of the coset representatives $\bc_1,\dots,\bc_d$ constructed in~\eqref{cosets} above has to be in $\Lambda \setminus  \left( \bigcup_{j=1}^m \Omega_j \right)$. Take $\bz$ to be that coset representative, then the bound on $|\bz|$ follows from~\eqref{cc_bnd}.
\endproof
\bigskip

\section{Angular distribution in~$\zed^n$}
\label{angle_sec}

In this section, we apply the results on small-norm zeros of integral quadratic equations to the problem of effective distribution of angles between vectors in~$\zed^n$. Specifically, let us write $\angle(\ba,\bb)$ for the angle between two vectors $\ba,\bb \in \zed^n$ and define
$$\Theta_n = \left\{ \angle(\ba,\bb) : \ba,\bb \in \zed^n \right\}$$
to be the set of all possible angles between such vectors. Fixing a particular vector $\ba \in \zed^n$, let us also write
$$\Theta_n(\ba) = \left\{ \angle(\ba,\bb) : \bb \in \zed^n \right\}.$$
It is established in~\cite{angles} that $\Theta_n = \Theta_n(\ba)$ for every $\ba \in \zed^n$ whenever $n = 2$ or $n \geq 4$, although this is not the case when $n=3$. In particular, when $n \geq 5$,
$$\Theta_n(\ba) = \{ \pi/2 \} \cup \left\{ \theta : \tan^2 \theta \in \que \right\}$$
for each $\ba \in \zed^n$. For $n \leq 4$, an integer vector making an assumed angle $\theta$ with $\ba$ can be fairly explicitly described, as shown in~\cite{angles}, so we will focus here on the situation $n \geq 5$ where the results of~\cite{angles} are less explicit (relying on Meyer's theorem instead; see, e.g., \S 6.1 of~\cite{cassels_qf}). Assuming $\theta \in \Theta_n(\ba)$, there can exist multiple vectors $\bb \in \zed^n$ such that $\angle(\ba,\bb) = \theta$. If $\theta = \pi/2$, then all such vectors are characterized as solutions to the linear equation $\sum_{i=1}^n a_ix_i = 0$, however the situation is more interesting when $\theta \neq \pi/2$. Here is our first observation in this direction.

\begin{lem} \label{angle_qf} Let $n \geq 5$ and $\ba = (a_1,\dots,a_n)^{\top} \in \zed^n$ be a nonzero vector and write $t = \|\ba\|^2 \in \zed$. Let $\pi/2 \neq \theta \in \Theta_n(\ba)$, so $\tan^2 \theta = q/p \in \que$. Then $\angle(\ba,\bb) = \theta$ for some $\bb \in \zed^n$ if and only if $\bb$ is a nontrivial zero of the quadratic form
\begin{equation}
\label{Qat}
Q_{\ba,\theta}(\bx) = pt \sum_{i=1}^n x_i^2 - (p+q) \left( \sum_{i=1}^n a_i x_i \right)^2.
\end{equation}
\end{lem} 

\proof
Notice that $\angle(\ba,\bb) = \theta$ if and only if
$$\cos \theta = \frac{\ba \cdot \bb}{\|\ba\| \|\bb\|},$$
where $\ba \cdot \bb = \sum_{i=1}^n a_i b_i$, $\|\ba\| = \sqrt{t}$, and $\|\bb\| = \left( \sum_{i=1}^n b_i^2 \right)^{1/2}$. Then
\begin{equation}
\label{pq_id}
\frac{q}{p} = \tan^2 \theta = \frac{1 - \cos^2 \theta}{\cos^2 \theta} = \frac{t \sum_{i=1}^n b_i^2 - \left( \sum_{i=1}^n a_i b_i \right)^2}{ \left( \sum_{i=1}^n a_i b_i \right)^2},
\end{equation}
which is equivalent to saying that $\bb$ is a zero of $Q_{\ba,\theta}(\bx)$.
\endproof

Now notice that~\eqref{Qat} can be expanded as 
$$Q_{\ba,\theta}(\bx) = pt \sum_{i=1}^n x_i^2 - (p+q) \sum_{i=1}^n \sum_{j=1}^n a_i a_j x_i x_j,$$
so
\begin{equation}
\label{Qat_norm}
|Q_{\ba,\theta}| \leq \max \{ pt, 2 (p+q) a_i a_j : 1 \leq i,j \leq n \} < 2 (p+q) \|\ba\|^2,
\end{equation}
and the symmetric coefficient matrix of $Q_{\ba,\theta}$ is
$$\A_{\ba,\theta} = \begin{pmatrix} p \|\ba\|^2 - (p+q) a_1^2 & -(p+q)a_1a_2 & \dots & -(p+q)a_1a_n \\
-(p+q)a_1a_2 & p \|\ba\|^2 - (p+q) a_2^2 & \dots & -(p+q)a_2a_n \\
\vdots & \vdots & \ddots & \vdots \\
-(p+q)a_1a_n & -(p+q)a_2a_n & \dots & p \|\ba\|^2 - (p+q) a_n^2 \end{pmatrix}.$$ 

\begin{lem} \label{A_det} $\det \A_{\ba,\theta} = -p^{n-1} q \|\ba\|^{2n}$.
\end{lem}

\proof
Notice that $\Theta_n(\ba) = \Theta_n$, i.e. it does not depend on $\ba \neq \bo$. This means that for any nonzero vector $\ba$ there exists $\bb$ so that~\eqref{pq_id} holds for the given $q/p$. Picking such $\ba$ and $\bb$, there exists a real orthogonal matrix $U_{\ba}$ such that
$$\ba = \|\ba\| U_{\ba} \be_1,$$
where $\be_1 = (1,0,\dots,0)$, and then we can write the vector $\bb = U_{\ba} \bb'$ for some appropriate $\bb' \in \zed^n$. Then $\theta = \angle(\ba,\bb) = \angle(\be_1,\bb')$, and we can rewrite~\eqref{pq_id} as follows:
\begin{eqnarray*}
\frac{q}{p} & = & \frac{\|\ba\|^2 \| \bb \|^2 - \left( \ba \cdot \bb \right)^2}{ \left( \ba \cdot \bb \right)^2} = \frac{\|\ba\|^2 \| U_{\ba} \bb' \|^2 -  \|\ba\|^2 \left( \be_1^{\top} (U_{\ba}^{\top} U_{\ba}) \bb' \right)^2}{ \|\ba\|^2 \left( \be_1^{\top} (U_{\ba}^{\top} U_{\ba}) \bb' \right)^2} \\ 
& = & \frac{\sum_{i=1}^n (b'_i)^2 - (b'_1)^2}{(b'_1)^2},
\end{eqnarray*}
so $\bb'$ is a zero of the quadratic form $Q_{\be_1,\theta}(\bx') = p \sum_{i=1}^n (x'_i)^2 - (p+q) (x'_1)^2$. In fact,
$$Q_{\ba,\theta}(\bx) = Q_{\be_1,\theta}(\|\ba\| U_{\ba} \bx') = \|\ba\|^2 Q_{\be_1,\theta}(U_{\ba} \bx'),$$
and so
$$\det \A_{\ba,\theta} = \|\ba\|^{2n} \det ( U_{\ba}^{\top} \A_{\be_1,\theta} U_{\ba}) = \|\ba\|^{2n} \det \A_{\be_1,\theta},$$
where
$$\A_{\be_1,\theta} = \begin{pmatrix} -q & 0 & \dots & 0 \\
0 & p & \dots & 0 \\
\vdots & \vdots & \ddots & \vdots \\
0 & 0 & \dots & p \end{pmatrix}.$$ 
The conclusion follows.
\endproof

This quadratic form is closely related to the nonsingular indefinite form defined in equation~(2.2) of~\cite{angles}. Further, by Lemma~\ref{A_det}, $\det \A_{\ba,\theta} \neq 0$, so $Q_{\ba,\theta}$ is also nonsingular, and it is isotropic over $\zed$ by Lemma~\ref{angle_qf} above, since we are assuming that there exists a vector $\bb \in \zed^n$ so that $\angle(\ba,\bb) = \theta$. Then we obtain the following effective observation as an immediate consequence of Cassels' bound~\eqref{cassels_bnd}.

\begin{cor} \label{angle_eff-1} Let $n \geq 5$ and $\bo \neq \ba \in \zed^n$. Let $\theta \in \Theta_n$, so that $\tan^2 \theta = q/p \in \que$. There exists a vector $\bb \in \zed^n$ such that $\angle(\ba,\bb) = \theta$ and
\begin{equation}
\label{ae_bnd}
|\bb| \ll \left( 2(p+q) \|\ba\|^2 \right)^{\frac{n-1}{2}},
\end{equation}
where the implied constant in the upper bound is explicitly computable and depends only on~$n$.
\end{cor}

\proof
Combining~\eqref{cassels_bnd} with~\eqref{Qat_norm} above, we obtain a point $\bb \in \zed^n$ with sup-norm bounded as in~\eqref{ae_bnd} so that $Q_{\ba,\theta}(\bb) = 0$. Hence, by Lemma~\ref{angle_qf}, $\angle(\ba,\bb) = \theta$.
\endproof

Further, we can use our Theorem~\ref{main} to obtain a bound on the sup-norm of a vector $\bb$ making a given angle with $\ba$ outside of a union of sublattices, provided such a vector exists.

\begin{cor} \label{angle_eff-2} Let $n \geq 5$ and $\bo \neq \ba \in \zed^n$. Let $\theta \in \Theta_n$, so that $\tan^2 \theta = q/p \in \que$. Let $\Omega_1,\dots,\Omega_m \subset \zed^n$ be proper sublattices of finite indices such that there exists $\bb \in \zed^n \setminus \left( \bigcup_{i=1}^m \Omega_i \right)$ with $\angle(\ba,\bb) = \theta$. Then there exists such a vector with
\begin{equation}
\label{ae_bnd-1}
|\bb| \ll \det (\Omega)^{2\rho(n)+1} \left( 2 (p+q) \|\ba\|^2 \right)^{\rho(n)},
\end{equation}
where $\Omega = \bigcap_{i=1}^m \Omega_i$ and the implied constant in the upper bound is explicitly computable and depends only on~$n$.
\end{cor}

\proof
Combining Theorem~\ref{main} with~\eqref{Qat_norm} above, we obtain a point $\bb \in \zed^n \setminus \left( \bigcup_{i=1}^m \Omega_i \right)$  with sup-norm bounded as in~\eqref{ae_bnd-1} so that $Q_{\ba,\theta}(\bb) = 0$. Hence, by Lemma~\ref{angle_qf}, $\angle(\ba,\bb) = \theta$.
\endproof

\bigskip

\bigskip

\bibliographystyle{plain}  

\begin{thebibliography}{10}


\bibitem{cassels_gn}
J. W. S. Cassels.
\newblock An Introduction to the Geometry of Numbers.
\newblock Springer-Verlag, 1959.

\bibitem{cassels:small}
J.~W.~S. Cassels.
\newblock Bounds for the least solutions of homogeneous quadratic equations.
\newblock {\em Proc. Cambridge Philos. Soc.}, 51:262--264, 1955.

\bibitem{cassels_qf}
J. W. S. Cassels.
\newblock Rational Quadratic Forms.
\newblock London Mathematical Society Monographs, 13. Academic Press, Inc., London-New York, 1978.

\bibitem{lf_chan}
W.~.K.~Chan and L.~Fukshansky.
\newblock Small representations of integers by integral quadratic forms.
\newblock {\em J. Number Theory}, 201 (2019), pp. 40--52.

\bibitem{dietmann}
R. Dietmann.
\newblock Small solutions of quadratic Diophantine equations.
\newblock {\em Proc. London Math. Soc.}, 86 (2003), no. 3, pp. 545--582.

\bibitem{lf1}
L.~Fukshansky.
\newblock Integral points of small height outside of a hypersurface.
\newblock {\em Monatsh. Math.}, 147 (2006), no. 1, pp. 25--41.

\bibitem{cassels_overview}
L.~Fukshansky.
\newblock Heights and quadratic forms: on {C}assels' theorem and its
  generalizations.
\newblock In W.~K. Chan, L.~Fukshansky, R.~Schulze-Pillot, and J.~D. Vaaler,
  editors, {\em Diophantine methods, lattices, and arithmetic theory of
  quadratic forms}, Contemp. Math., 587, pages 77--94. Amer. Math. Soc.,
  Providence, RI, 2013.

\bibitem{gruber_lek}
P. M. Gruber and C. G. Lekkerkerker.
\newblock Geometry of Numbers.
\newblock North-Holland Publishing Co., 1987.

\bibitem{henk_thiel}
M. Henk and C. Thiel.
\newblock Restricted successive minima.
\newblock {\em Pacific J. Math.}, 269 (2014), no. 2, pp. 341--354.

\bibitem{kornhauser}
D. M. Kornhauser.
\newblock On the smallest solution to the general binary quadratic diophantine equation.
\newblock {\em Acta Arith.}, 55 (1990), no. 1, pp. 83--94.

\bibitem{martinet}
J. Martinet.
\newblock Perfect Lattices in Euclidean Spaces.
\newblock Springer-Verlag, 2003.

\bibitem{masser}
D.~W. Masser.
\newblock How to solve a quadratic equation in rationals.
\newblock {\em Bull. London Math. Soc.}, 30(1):24--28, 1998.

\bibitem{angles}
K.~Yamamoto.
\newblock On exceptionality of dimension three in terms of lattice angles.
\newblock {\em arXiv:2304.02299}, 2023.



\end{thebibliography}

\end{document}